# Structural System Identification via Validation and Adaptation

Cristian López[1], Keegan J. Moore[1,2]*

[1]*Department of Mechanical and Materials Engineering, University of Nebraska-Lincoln, Lincoln, NE 68588, United States*
[2]*Daniel Guggenheim School of Aerospace Engineering, Georgia Institute of Technology, Atlanta, GA 30332, United States*

*Corresponding author: (K.J. Moore)
E-mail address: kmoore@gatech.edu

## ABSTRACT

Estimating the governing equation parameter values is essential for integrating experimental data with scientific theory to understand, validate, and predict the dynamics of complex systems. In this work, we propose a new method for structural system identification (SI), uncertainty quantification, and validation directly from data. Inspired by generative modeling frameworks, a neural network maps random noise to physically meaningful parameters. These parameters are then used in the known equation of motion to obtain fake accelerations, which are compared to real training data via a mean square error loss. To simultaneously validate the learned parameters, we use independent validation datasets. The generated accelerations from these datasets are evaluated by a discriminator network, which determines whether the output is real or fake, and guides the parameter-generator network. Analytical and real experiments show the parameter estimation accuracy and model validation for different nonlinear structural systems.

**Keywords:** data-driven, generative adversarial networks, parameter estimation, uncertainty quantification, vibrating structures.

## 1. Introduction

Understanding and managing vibrations is essential for ensuring the longevity and efficiency of structures and devices [1, 2]. System identification, an important methodology in structural dynamics, uses, for example, vibration data to derive mathematical models that capture the system behavior, enable parameter estimation (e.g., mass, stiffness, and damping), and predict dynamic response [3, 4]. This approach is fundamental in various fields for enhancing safety and performance in applications such as structural health monitoring, where real-time sensing systems detect damage and assess the load-carrying capacity of infrastructure [5], in vibration control systems, as accurate models can mitigate adverse effects of external forces, ensuring structural integrity and stability during environmental loads like seismic activity or wind [6], in noise control, such that precise modeling of a system's behavior enables the design of effective noise reduction strategies [7], etc.

Traditionally, approaches to vibration analysis have relied on theoretical models validated through experimental testing designed around those models [8]. Methods based on linear assumptions and Fourier techniques often inadequately capture the nonlinearity and nonstationary dynamics of the systems [9]. Modern data-driven methodologies now can characterize the structural dynamics directly from measurements [10], offering new opportunities to understand these complex behaviors [11, 12].

Kerschen et al. [3] and Noël et. al [4] reviewed advances in nonlinear SI in structural systems. Methods frequently applied for nonlinear SI can be divided into three types: *parametric methods* require a predefined mathematical model, with the identification process focusing on estimating its coefficients. Commonly, these methods include moving average models [13, 14], Kalman filter [15], Bayesian methods [16, 17], polynomial nonlinear state-space system models [18], using time-series based [19, 20], and moving horizon optimization [21] methods, infusing physics into neural networks (NNs) [22–26], etc.; *non-parametric methods* do not rely upon prior knowledge of the governing dynamics, instead, they infer the system's equations of motion from data. Numerous methods have been proposed, e.g., the restoring force surface method [27, 28], employing genetic programming [29–31], NNs [32–36], using a slow-flow model of the dynamics [37, 38], a Nonlinear Identification through eXtended Outputs method [39], employing symbolic regression and genetic programming [40], a data-driven approach that integrates machine learning and symbolic regression [41], etc.; *semi-parametric methods* integrate prior but incomplete knowledge of the system with data-driven methods. In this category, several methods have been proposed, such as the piecewise-linear RFS method [42], the sparse identification of nonlinear dynamics (SINDy) method [43] applied to systems with inelastic/hysteresis phenomena [12], a Bayesian framework embedded in SINDy [44], the characteristic nonlinear SI method [45], a Hamiltonian-constrained autoencoder [46], an energy-based approach [47, 48].

Over the past few years, as artificial intelligence has advanced, data-driven methods have become widely adopted for SI tasks. Cunha et al. [49] summarized recent computational intelligence techniques for nonlinear dynamic SI, and Quaranta et al [50] reviewed the application of machine learning approaches to structural dynamics and vibroacoustic. These advances exploit the capability of machine learning frameworks to reveal intricate patterns and dependencies within data, eliminating the necessity for explicitly defined mathematical formulations. For instance, a feed-forward NN was employed to estimate the restoring force of a Duffing oscillator [32], a wavelet NN model was developed for structural system identification [33], a recurrent NN was applied to determine the impact forces of nonlinear structures [34], a convolutional NN was utilized to predict structural dynamic response and perform system identification [35], a symbolic NN was proposed to derive mathematical models for the governing dynamics [36]. In many cases, SI benefits from incorporating available physics knowledge, where some parameters or governing equations are (partially) known. Embedding this knowledge into the training process allows the model to balance data-driven insights with physical consistency and constraints. For example, [22] utilized a long short-term memory network to predict the nonlinear structural responses of systems subjected to ground motion excitations, while [23] used this type of network for nonlinear hysteretic parameter identification; [24] combined the Runge-Kutta integration scheme with physics-informed NNs to identify parameters and model nonlinear dynamical systems, [25] employed generative adversarial networks (GANs) for structural parameter identification, [26] used fully connected networks to estimate nonlinear structural responses, [46] utilized physical constraints from Hamiltonian mechanics in an autoencoder to discover interpretable models of structural dynamics, etc. Among these architectures, GANs have emerged as a powerful technique that can generate realistic and diverse data [51]. It consists of two neural networks, a generator and a discriminator, that compete during training; the generator learns to produce realistic synthetic data, while the discriminator learns to distinguish between real and generated data, with both networks improving through this adversarial process. This technique is of interest for applications such as image generation [52], text generation [53, 54] and music generation [53], mechanical vibration time series data augmentation [55], visual anomaly discovery [56], superconductors with high transition temperatures [57], to solve stochastic differential equations [58], mechanical systems fault detection [59], etc. In the field of structural mechanics, GANs have been employed to determine correlations between vibration measurements and the parameters of nonlinear models, aiming for model updating [60], identifying and modeling parameterized nonlinear systems [61], nonlinear modal analysis [62], and parameter identification [25].

In this paper, inspired by GANs, we introduce a new parametric approach for structural SI coined the System Identification via Validation and Adaptation (SIVA) method. This technique computes the system's parameters with training acceleration time series while simultaneously validating them with acceleration time series not used in the parameter estimation procedure. This is feasible because the calculated parameters are used to generate predicted acceleration time series, which are then compared with the actual training accelerations to update the parameters iteratively. The estimated parameters are then leveraged to generate acceleration time series from the validation data set and then sent to a discriminator, which is used to perform validation. As such, the SIVA method estimates and validates the differential equation's parameters directly from measured responses of vibrating structures. Additionally, once convergence has been reached, uncertainty quantification is performed on the estimated parameters. Moreover, we will demonstrate that the method is robust to overtraining. The method's performance is shown by utilizing both nonlinear analytical and experimentally measured systems.

## 2. The Proposed Method

A similar strategy was presented in [25], where they combined learning capabilities of GANs with physics-based constraints for structural parameter identification. Specifically, they incorporated physical knowledge by using generators for structural parameters and displacements, and by formulating physics-based loss functions that link the data generator with parameter generators through the governing equations of motion. Our work is different from this approach because: 1) we use an architecture that identifies the system parameters (damping, stiffness) and use them in a physical model to produce accelerations, 2) the discriminator tries to distinguish between real and generated accelerations, and 3) we perform validation using different datasets than the identification data. Details of the proposed method are presented in the following discussion with a flowchart for the method shown in Fig. 1. Throughout this paper, we use bold lowercase letters for vectors and bold capital letters for matrices.

## 2.1. Procedure of the proposed methodology

Our objective is to simultaneously learn, validate, and quantify uncertainty for a representative mathematical model using measured data and adversarial processes based on GANs. To achieve this, we employ the following methodology:

1. We begin by processing the data and assuming a mathematical model $f(q, \dot{q}; \lambda)$, where $q$ and $\dot{q}$ represent the displacement and velocity vectors, respectively, and $\lambda$ represents the unknown parameters of the model.

2. To identify the parameters of the model, a parameter-generator network ($P$) transforms a batch of independent random noise $z$ into meaningful physical parameter values, which are then incorporated into $f$.

3. The optimization is guided by both the mean square error (MSE) between the real accelerations $\ddot{q}_{tr}$ and model-generated accelerations $\tilde{\ddot{q}}_{tr} = f(\dot{q}_{tr}, q_{tr})$; and the adversarial loss, which evaluates how well the generated parameters, and therefore the generated accelerations, can 'fool' the discriminator network ($D$), such that $\tilde{\ddot{q}}_{val}$ are classified as real. This loss provides direct feedback to the parameter generator network using the identification dataset.

4. To simultaneously validate the model, we introduce independent validation data collected under different experimental conditions from the training dataset. The validation data ($\dot{q}_{val}, q_{val}$) is feed into the identified model, and the resulting output $\tilde{\ddot{q}}_{val}$ is evaluated by a discriminator network, which determines whether the output is real or fake. This classification is used to update the discriminator and the parameter-generator network.

5. Once convergence is achieved, uncertainty quantification (UQ) is performed on the identified parameters, which can be performed either by continuing the training process and harvesting the parameter values from additional epochs or by applying the parameter generator to new noise data without training. The first approach is possible due to the robustness to overtraining that the methodology exhibits. The second approach allows one to halt the training immediately after convergence and perform UQ without the cost associated with the training process.

Through the adversarial process, both networks progressively improve, ultimately providing a framework with inherent validation. While tracking the convergence of the parameters, if computational resources are available, the optimal parameter combination can be obtained by simulating the displacements, $\tilde{q}_{tr}$, of the system and comparing them with the real displacements, $q_{tr}$. The selected combination is based on minimizing the MSE between these two displacement vectors, and can be done using identification data, validation data, or both.

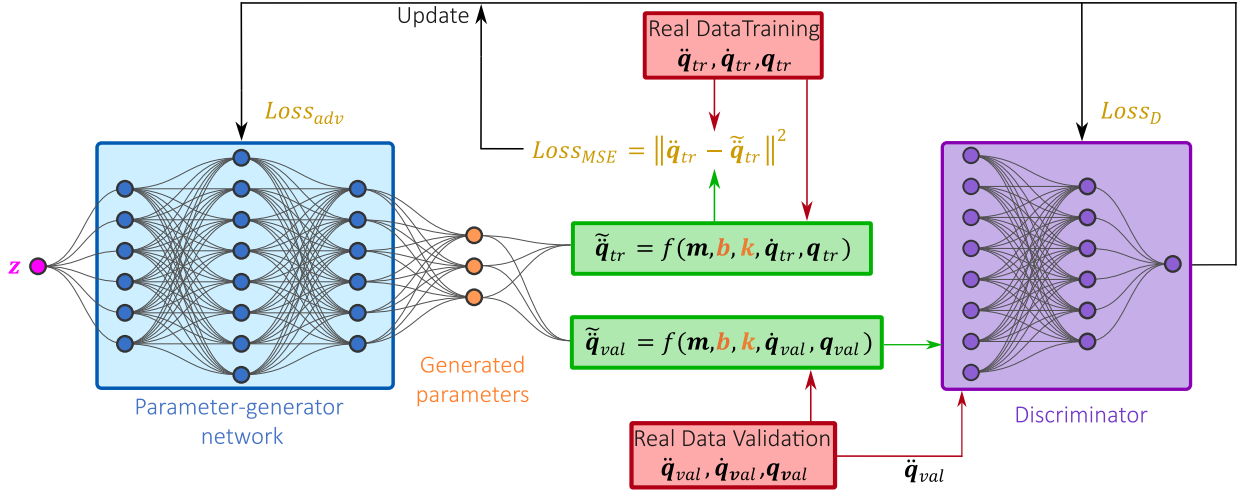

**Fig. 1.** Flowchart of the proposed SIVA method for parametric system identification.

## 2.2. Data processing

In this work, we assume that the masses are measured and that enough of the physics is known to be able to propose a reduced-order mathematical model for the system, but the parameters of the model are unknown and need to be identified from measured response data. The assumption of known masses can be relaxed and mass-scaled parameter can be identified instead. For experimental measurements obtained from accelerometers, displacements and velocities can be calculated by applying numerical integration and filtering, following the method in [63]. If the signals are embedded in noise, the proposed approach may not be able to accurately identify a model for the system or generate meaningful acceleration time series. However, this was not explicitly tested in this research and is left open for future research. Furthermore, in this work, experimental measurements are bandpass filtered to minimize the effects of noise at low and high frequencies.

## 2.3. Network architecture

The proposed framework contains two neural networks adapted for structural parameter identification. The parameter-generator network transforms random input noise $z$ into physically meaningful system parameters $\lambda = [b, k]$, where the parameters have been divided into damping and stiffness coefficients. This network has a deep, fully connected architecture with four layers ($64 \rightarrow 32 \rightarrow 16 \rightarrow n$ neurons, with $n$ as the total number of parameters to identify), each followed by `LeakyReLU` activation functions with a negative slope 0.2, except for the final layer, which uses a linear activation function. A linear activation function is used at the end to preserve negative values if needed. For systems where the parameters span several orders of magnitude, we apply a scientific notation transformation during learning. Specifically, a physical parameter used by the model is $k = a \cdot 10^b$, where both a and b are values to learn for each desired parameter. This approach acts as a form of normalization that enhances the effectiveness of the parameter generator.

The discriminator network, $D$, is similar to the parameter-generator network except that it has three fully connected layers ($64 \rightarrow 32 \rightarrow 1$ neurons) and applies the `Sigmoid` function to the output layer. The proposed approach was developed in Python 3.12.7 with Pytorch 2.6.0 platform and also tested in Python 3.10.16 using Pytorch 2.6.0. The results in this research were produced with the former setup. We use Adam as the optimizer with a learning rate of $10^{-4}$ with a batch size of 300.

## 2.4. Physics-based modeling

For a multiple-degree-of-freedom system, from Newton's second law, we compute the accelerations $\ddot{q}$ based on the known mass matrix $M$, calculated parameters, and measured state variables ($q$ and $\dot{q}$) as

$$\ddot{q} = M^{-1}[-B(q,\dot{q}) - K(q) + F(t)], \tag{1}$$

where the matrices $B$, $K$, and $F$ represent internal non-conservative, internal conservative, and external applied forces, respectively. Note that Eq. (1) can also be applied to single-degree-of-freedom systems. Once the system's parameters have been obtained from the parameter-generator network we obtain the fake accelerations $\widetilde{\ddot{q}}$ using Eq. (1), for both the training and validation datasets.

### 2.4.1. Loss functions

The networks are trained following conventional adversarial learning [51]. The discriminator is trained to maximize the probability of correctly classifying real accelerations as real and generated accelerations as fake. Simultaneously, the parameter-generator is trained to minimize the probability that the discriminator accurately classifies the generated data. In this work, in addition to the adversarial learning, we include the MSE between generated and real accelerations.

The discriminator loss employs the traditional binary cross-entropy criterion:

$$\mathcal{L}_D = -\mathbb{E}_{\ddot{q} \sim p(\ddot{q})}[\log(D(\ddot{q}))] - \mathbb{E}_{z \sim p(z)}\left[\log\left(1 - D(\widetilde{\ddot{q}})\right)\right], \tag{2}$$

where $z$ is the latent variable sampled from the prior distribution $p_z(z)$. We use a normal distribution $\mathcal{N}(0,1)$, with the `randn` function, and seed 42 for reproducibility. $D(\ddot{q})$ represents the probability that $\ddot{q}$ came from the real acceleration data (samples from the real distribution $p(\ddot{q})$), and $D(\widetilde{\ddot{q}})$ denotes the probability that $\widetilde{\ddot{q}}$ came from generated data, which is obtained as $f(q, \dot{q}; m, b, k)$, with $b, k = P(z)$. It follows the two-step training process described in the GAN's literature [51]. We train the discriminator network to distinguish between real accelerations from validation data and fake accelerations generated by the model, $f$, using the standard GAN approach where real samples are labeled as 1 and generated samples as 0 [64].

The optimal discriminator loss is obtained when the discriminator can no longer distinguish between real and generated samples as

$$\mathcal{L}_D = -\mathbb{E}_{\ddot{q} \sim p(\ddot{q})}[\log(0.5)] - \mathbb{E}_{z \sim p(z)}[\log(1 - 0.5)] = \log(4) = 1.386. \tag{3}$$

The parameter-generator loss combines two terms:

$$\mathcal{L}_P = \mathcal{L}_{adv} + \gamma \mathcal{L}_{\text{MSE}}, \tag{4a}$$

$$\mathcal{L}_P = -\mathbb{E}_{z \sim p(z)}\left[\log\left(D(\widetilde{\tilde{q}})\right)\right] + \gamma \mathbb{E}_{\ddot{q} \sim p(\ddot{q})}\left[\|\ddot{q} - \widetilde{\tilde{q}}\|^2\right], \tag{4b}$$

where the first term is the adversarial loss, encouraging the parameter-generator to provide values such that the generator block produces realistic examples, while the second term ensures that the generated acceleration samples are similar to the real ones. The hyperparameter $\gamma$ balances the importance of these two terms. We train the parameter-generator while keeping the discriminator's parameters fixed. For the first term in Eq. (4), the target labels for the discriminator's output are set to 1s, meaning the parameter-generator is trained to transform noise vectors into meaningful parameters, that are used in the model, such that the discriminator accepts $\widetilde{\tilde{q}}$ as real.

### 2.4.2. Uncertainty quantification

To assess the reliability of the estimated parameters, we perform UQ using the parameter values obtained during model training. Specifically, after convergence is reached, we store the estimated parameters at each epoch, though this could also be done at the batch level for more detailed insights if desired. We then analyze the distribution of each parameter by fitting a normal distribution to the obtained parameter samples using MATALB's `fitdist` function (this could be done in Python as well using SciPy if desired). After that, we evaluate the probability density function using the `pdf` function in MATLAB (with evenly spaced points $\pm 6$ standard deviations from the mean), which computes the probability density values based on the fitted distribution parameters. This function provides a smooth approximation of the probability density function for each parameter, that allows us to evaluate the spread and concentration of parameter estimates over training.

This work is inspired by the standard GAN by incorporating physics-based inputs, specifically the governing equation of motion for the chosen model, and the measured displacements and velocities to generate the corresponding accelerations. This approach uses the adversarial framework to learn the physical parameters that govern the system dynamics. Unlike most SI methods, the proposed approach not only performs data-driven and physics-based parameter identification but also UQ and model validation simultaneously.

## 3. Method Demonstration

### 3.1. Analytical single-degree-of-freedom nonlinear oscillator.

To demonstrate SIVA, we utilize data simulated from a Duffing oscillator, presented in [47], with the equation of motion

$$m\ddot{x} + b\dot{x} + b_{nl}x^2\dot{x} + kx + k_{nl}x^3 = 0, \quad x(0) = 0 \text{ m}, \dot{x}(0) = 5 \text{ m/s}, \tag{5}$$

where $m = 0.05$ kg is the mass, $b = 0.5$ Ns/m is the linear damping coefficient, $b_{nl} = 4000$ Ns/m³ is the nonlinear damping coefficient, $k = 300$ N/m is the linear stiffness, and $k_{nl} = 3 \times 10^8$ N/m³ is the nonlinear stiffness, and $x = q$ is the displacement coordinate. To produce the simulated time series, we numerically solve Eq. (5) using MATLAB's `ode45` with relative and absolute tolerances set to $10^{-8}$. We set the time duration to 1 s and the sampling rate to 10 kHz ($dt = 0.0001$). Figure 2 displays the simulated time series along with its continuous wavelet transform (CWT) spectrum (normalized such that the maximum amplitude is 1) [65, 66] and Fourier spectrum. For the validation time series, we use $x(0) = 0$ m, $\dot{x}(0) = 8$ and 3 m/s, as shown also in Fig. 2.

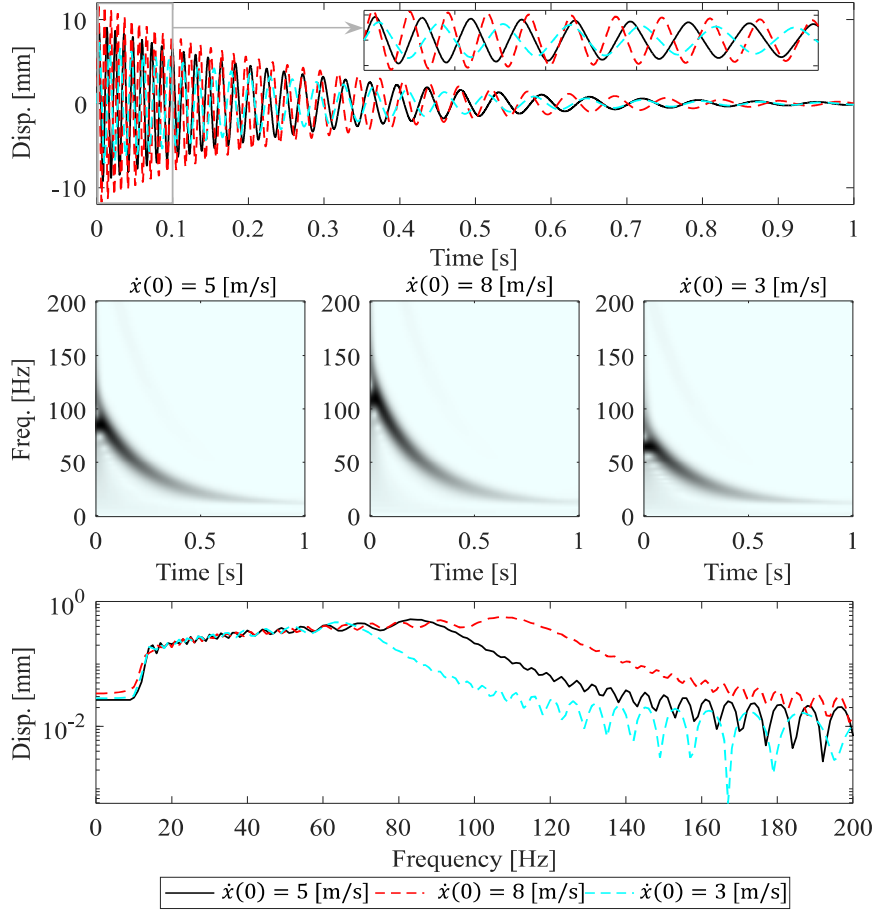

**Fig. 2.** Simulated displacement response of the Duffing oscillator.

Figure 3(a) displays the training process of the proposed methodology, where significant fluctuations in the loss values can be seen during the initial epochs as the parameter-generator and discriminator compete to improve. After around 400 epochs, the training converges, with the losses attempting to stabilize. On the left panel, the discriminator loss approaches the expected theoretical value of $\ln(4)$. Although the discriminator loss has not exactly converged to the theoretical value, its value is sufficiently close, and the parameter values also show sufficient convergence. This indicates that the parameter-extractor has learned to produce precise parameters that allow the model to obtain accelerations capable of fooling the discriminator. We note that if the training process is carried on, the mean values of the parameters for epoch 2000 are $b = 0.50174$ Ns/m, $b_{nl} = 4152.8$ Ns/m³, $k = 298.29$ N/m, and $k_{nl} = 2.9727 \times 10^8$ N/m³, which differ only slightly from the mean parameters identified from epoch 1000, which are $b = 0.48930$ Ns/m, $b_{nl} = 4058.3$ Ns/m³, $k = 299.98$ N/m, and $k_{nl} = 2.9198 \times 10^8$ N/m³. While on the right panel, the first term of Eq. (4) converges to approximately $\ln(4)/2$, and the MSE decreases to a relatively small value, indicating a good agreement between the simulated and generated acceleration signals.

In the following, we consider three approaches, ordered from least to greatest computational cost, for determining the model parameters:

- Approach I: After training, each parameter is sampled $N$ times [25] from the parameter-generator (we use $N = 1000$ in this work, but any value $N \geq 1$ can be used). The final model parameters are then obtained by computing the mean of the $N$ samples. We present these values in Table 1.

- Approach II: After convergence is reached, the training process is continued and the parameter values from each epoch stored until the total number of desired epochs is achieved. The mean value is computed for each parameter and used for the final model parameters. We applied this approach for the parameters identified after epoch 400 and provided those values in Table 1. This approach is included to demonstrate the robustness of the methodology to overtraining. However, in practice, it is much more economical to use Approach I unless the methodology is applied continuously and ingests new data as it becomes available.

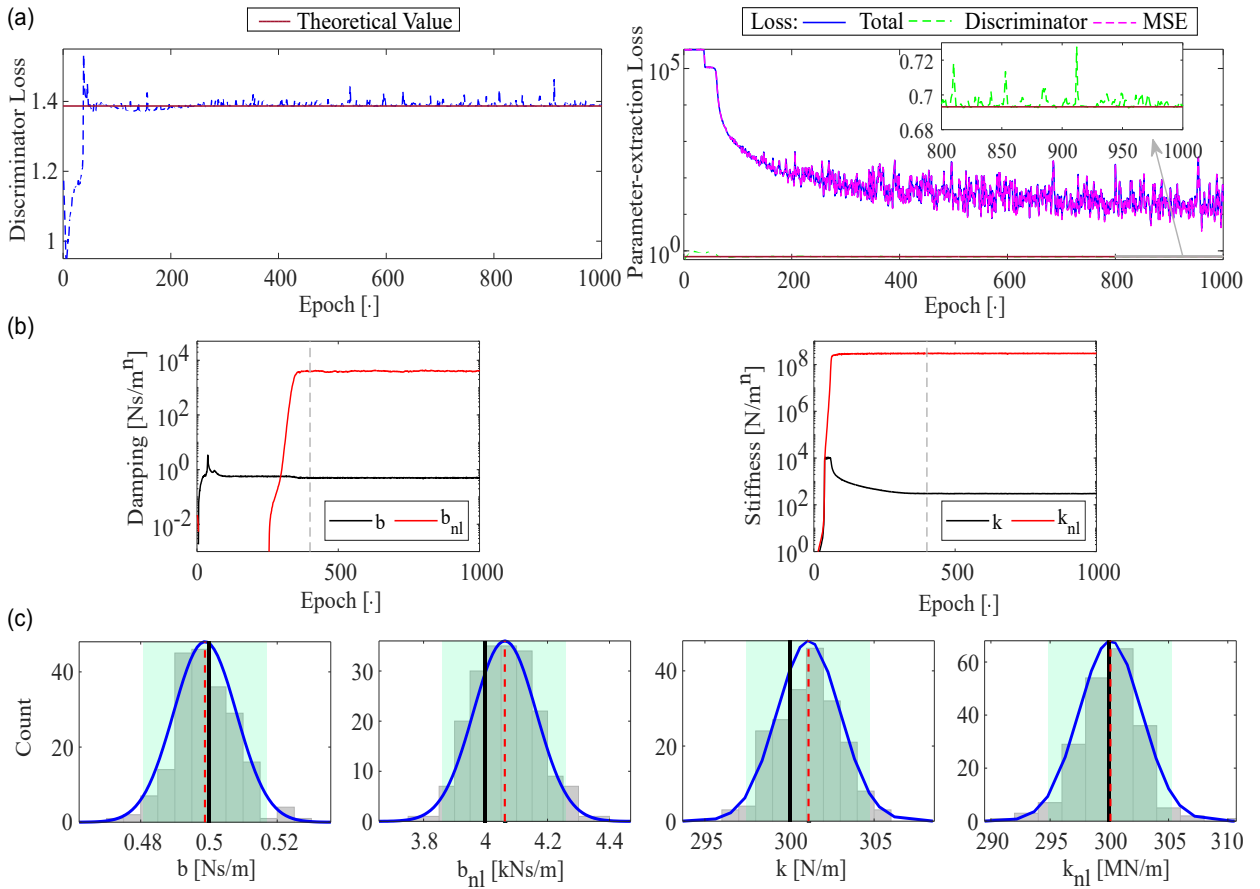

**Fig. 3.** Training dynamics of the proposed SIVA: (a) Losses, (b) Identified parameters distributions of the Duffing Oscillator system.

- Approach III: During training, if computational cost is not an issue, one can incorporate an explicit time integration and time series comparison to obtain the optimal set of parameters from the identified set. We simulated the response of the system using the full parameter set produced using Approach II (epochs 400 to 1000) by solving an initial value problem through numerical integration for each combination of parameters. We note that a time-series optimization routine, such as those used in [19, 20, 67, 68], or a more sophisticated sampling approach (e.g., Latin hypercube sampling) could be used instead of simulating every parameter combination. The resulting simulated displacements ($\tilde{x}$) are compared to the true system response ($x$) using MSE as the metric (MSE = $\frac{1}{J}\sum_{j=1}^{J}(x_j - \tilde{x}_j)^2$). We maintain and update the best parameter set when a lower MSE is achieved, ensuring that the identified parameters satisfy the adversarial training criteria and accurately reproduce the system's physical behavior as

$$\boldsymbol{b}, \boldsymbol{k} = \arg\min_{\boldsymbol{b},\boldsymbol{k}} MSE(\boldsymbol{b}, \boldsymbol{k}), \qquad (6)$$

where $\boldsymbol{b}, \boldsymbol{k}$ are the optimal damping and stiffness parameters. We determined the optimal parameter values for Eq. (6), identified by the SIVA method and presented in Table 1.

We note that in addition to these approaches, one could use the parameters that produced the lowest MSE in training for the final model parameters as well if desired. In our testing, we found that these parameters differed than those produced by Approach III. Furthermore, any of these approaches could be used to produce initial guesses for optimization routines that update the parameters to reproduce the measured time series, such as in [19, 20, 67, 68]. Lastly, UQ can be performed for both Approaches I and II by fitting a distribution to sampled set of parameter values. However, Approach III does not provide a clear path for UQ, since it reduces the full set of values down to a single combination.

To simulate the system, we use the `solve_ivp` function in SciPy with relative and absolute tolerances set to $10^{-8}$. Alternatively, MATLAB® can be used to simulate the system response, and we observed faster run times using `ode45` in MATLAB over `solve_ivp`. Solving Eq. (5), we found that `ode45` took 0.09651 s on average while `solve_ivp` took 0.2124 s on average, such that MATLAB® is approximately 2.2 times faster than SciPy. These tolerance values are chosen because they reduce numerical errors without considerable computational cost on modern computers, but higher tolerance values can be used if computational cost is a limiting factor. Thus, to visualize the variability and accuracy relative to the true system, from the collection of simulated signals, we compute a mean response along with a shaded region representing the 95% confidence interval (CI) around the mean.

To demonstrate UQ with the proposed methodology, we fit a normal distribution to the parameters obtained using Approach II and note that this could also be done with the parameters from Approach I. From the obtained distributions in blue of Figure 3(c), we can see that the estimated parameters are concentrated around their respective mean values with relatively small standard deviations, indicating a consistent estimation for each parameter. In each subplot, the dashed red lines indicate the mean value, the green-shaded region represents the 95% confidence interval (green), and the exact value is shown as the solid black line.

**Table 1.** Comparison of the coefficients of the one-degree-of-freedom model.

| | | SINDy | EDDI | SIVA | | |
| --- | --- | --- | --- | --- | --- | --- |
| | | | | Approach I | Approach II | Approach III |
| Coefficient | Exact | Identified | Identified | Identified | Identified | Identified |
| $b$ [Ns/m] | 0.5 | 0.50089 | 0.50001 | 0.48832 | 0.49961 | 0.49999 |
| $b_{nl}$ [Ns/m$^3$] | 4000 | 3946.1 | 3998.8 | 4052.5 | 4004.4 | 3996.4 |
| $k$ [N/m] | 300 | 316.27 | 315.06 | 299.9 | 301.84 | 298.64 |
| $k_{nl}$ [N/m$^3$] | $3 \times 10^8$ | $2.9934 \times 10^8$ | $2.9937 \times 10^8$ | $2.9142 \times 10^8$ | $2.9997 \times 10^8$ | $2.9999 \times 10^8$ |
| MSE | – | $9.235 \times 10^{-3}$ | 0.01241 | 0.2835 | $3.262 \times 10^{-3}$ | $1.100 \times 10^{-4}$ |

For comparison to state-of-the-art, we also report the values obtained using the SINDy [43] and EDDI [47] methods for the terms present in Eq. (5) as the candidate functions. To quantify the success of the identification, we calculate the MSE between the real displacement, $x$, and the simulated time series signal, $\tilde{x}$. While EDDI performs well in identifying the parameters, SINDy proves to be more accurate than EDDI in this example. As observed, the proposed approach – particularly in Approach III – provides the best match for the displacement time series, but this result comes at the cost of increased computational time. The calculation time employed by the proposed SIVA for Approach II for 1000 epochs is 8 min 57 s. If we implement Approach III, this time is extended to 11 min 48 s. The calculation times for SINDy and EDDI are a fraction of a second; however, neither method incorporates validation or uncertainty quantification. Furthermore, the proposed method easily adapts to include additional datasets for either identification or validation cases. The training of SIVA was performed on a PC with an Intel Core i7-8700 CPU @ 3.2 GHz, 32 GB RAM, Windows 10, and 64-bit operating system.

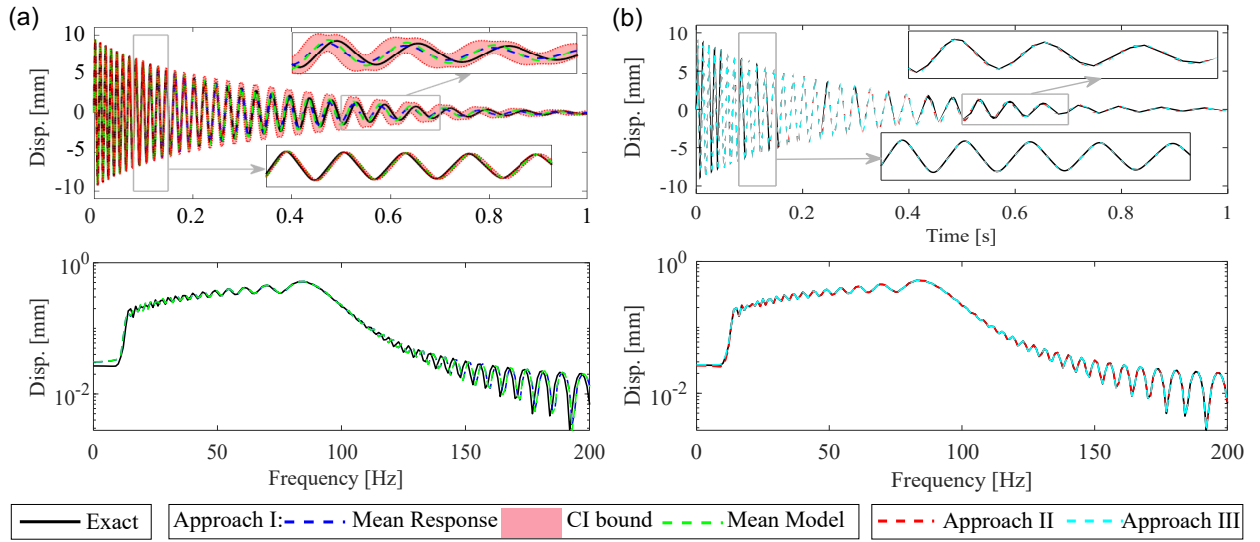

**Fig. 4.** Comparison of the displacement responses and Fourier spectra for the exact system and the identified model for $x(0) = 0$ m, $\dot{x}(0) = 5$ m/s. Simulation of signals from parameters obtained: (a) Approach I and (b) Approaches II and III.

Figure 4(a) reports both results from Approach I. In blue, we show the mean response signal generated from 1000 simulated signals obtained from the parameters provided by the trained parameter-generator network, along with the 95% confidence (pink); while in green, we present the simulated signal after using the mean values of the 1000 parameter sets. From the comparison of the displacement time series and Fourier spectra for the exact system and the responses from Approach I, we can see that the nonlinear part is well-captured, and the narrow confidence bounds indicate low variability for the sampled parameters. However, the response for both models deviates from the exact response at later times, which is caused by the mismatch in the damping. Using the identified parameters from Approaches II and III from Table 1, we simulated the response of the system by applying the same initial conditions as the exact system for training. As shown in Fig. 4(b), the identified model, particularly Approach III, closely matches the response of the exact system, confirming the capability of the proposed method to identify the system's parameters. Given the relative simplicity of the system and the absence of noise or external influences, a strong agreement is to be expected. Despite that, the results show the capabilities of the proposed method in determining the parameters of strongly nonlinear systems. Additionally, we expected a close match because the exact signal was part of the identification process. To further test the identified parameters of Approach III, we use the validation dataset. As depicted in Fig. 5, there is strong agreement between the exact signals and those produced by SIVA.

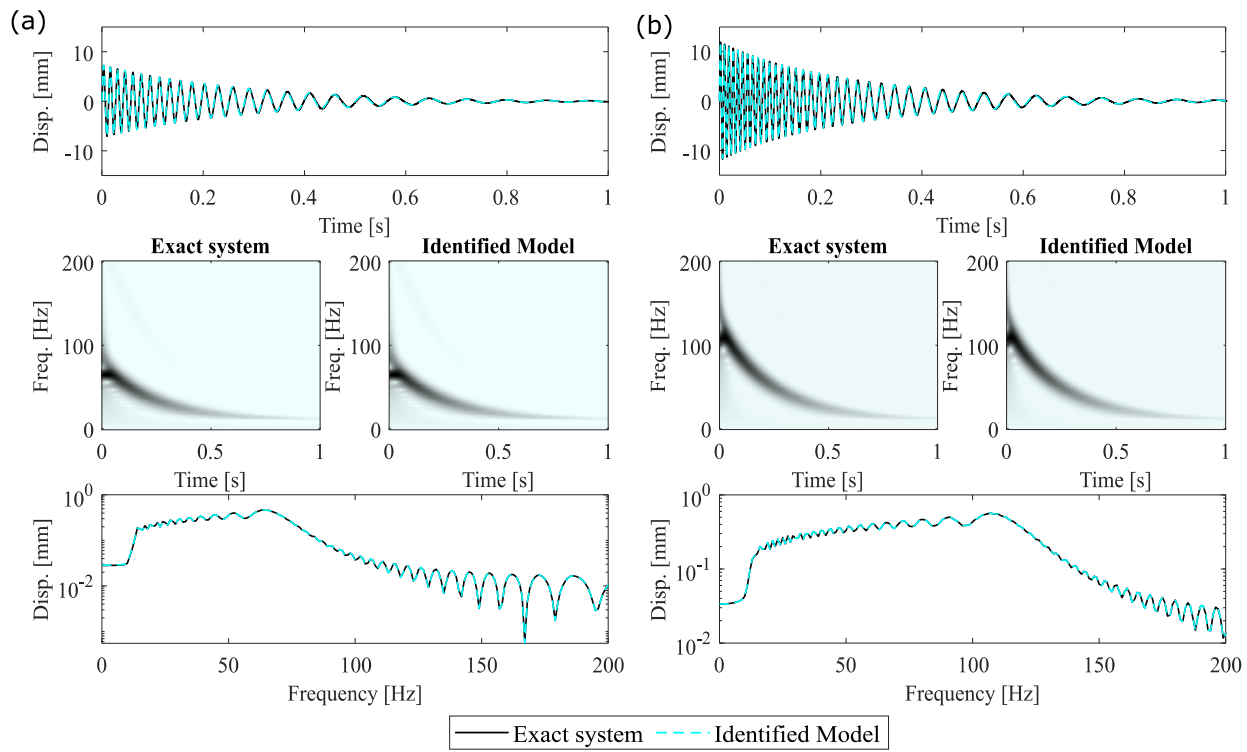

**Fig. 5.** Comparison of the exact time series and predicted responses for initial conditions $x(0) = 0$ m, (a) $\dot{x}(0) = 3$ m/s and (b) $\dot{x}(0) = 8$ m/s, used for validation, employing Approach III.

## 3.2. Experimental nonlinear coupled oscillators

The SIVA method is applied to a laboratory system, previously studied in [67, 68], to demonstrate its effectiveness when using real measured vibration signals. The experimental setup (Fig. 6) consists of two main components: the linear oscillator (LO) and the nonlinear oscillator (NO), both made of aluminum. The LO is grounded to an optical table using steel flexures and L-shaped brackets with the flexures acting as linear grounding springs. The mass of the assembled LO including coupling flexures is 1.37 kg. The NO is coupled to the LO using thin steel flexures and wires, which are clamped to aluminum anchors on the LO and NO. The steel flexures introduce weak linear stiffness, while the wires introduce strong stiffness nonlinearity due to geometric effects as well as additional weak linear stiffness. To balance the mass between the LO and NO, additional aluminum plates and washers were added to the NO, resulting in a mass equal to that of the LO. For detailed specifications of components, materials, and dimensions, please see [67, 68]. In reference [67], the experimental measurements were conducted in two phases: (1) the response of the LO alone was measured without the NO installed to identify a linear model; and (2) the free response of the coupled system (LO-NO) was measured for two forcing scenarios, where either the LO or NO was excited. Excitation was applied using a PCB Piezotronics modal impact hammer, and the free response was recorded using PCB accelerometers at a sampling rate of 2048 Hz for 8 s using Data Physics Abacus 906 hardware and software. The acceleration time series were band-pass filtered between 4 and 50 Hz using a third-order Butterworth filter implemented in MATLAB's `filtfilt` function, which avoids phase shift by applying the filter forwards then backwards. The velocities are obtained by numerically integrating the measured accelerations and then applying the third-order Butterworth high-pass filter with a cutoff frequency of 3 Hz. The displacements are computed using the same procedure but using the velocities instead.

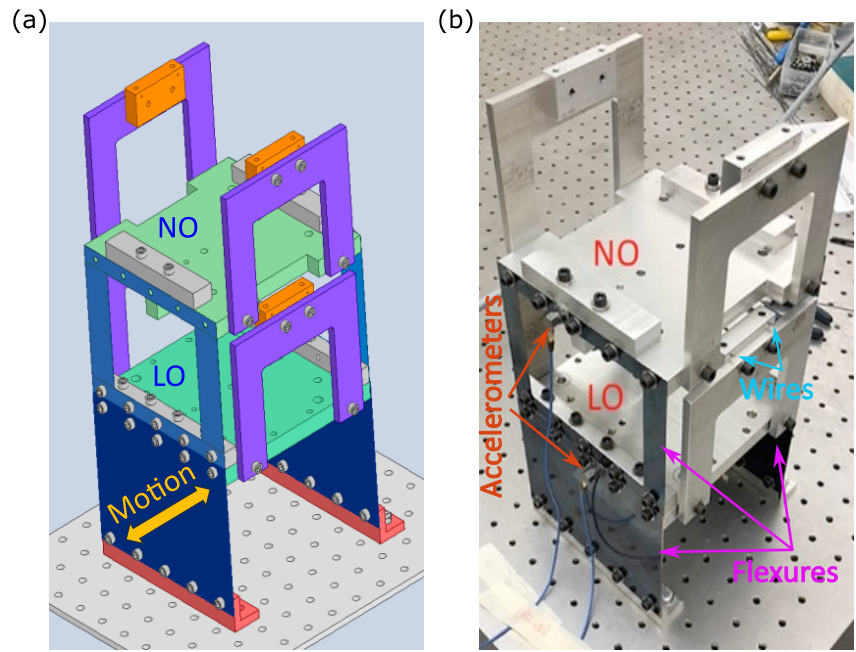

**Fig. 6**. (a) CAD model, (b) Instrumented two-degree-of-freedom system with flexures, wires, masses, and C-shaped brackets.

Figures 7(a), (c), and (e) show the forces applied to the system, whereas Figs. 7(b), (d), and (f) present the displacement responses and the CWT spectra for each applied force, respectively. The resulting response is nonlinear as seen in the early times and the system decays to linear

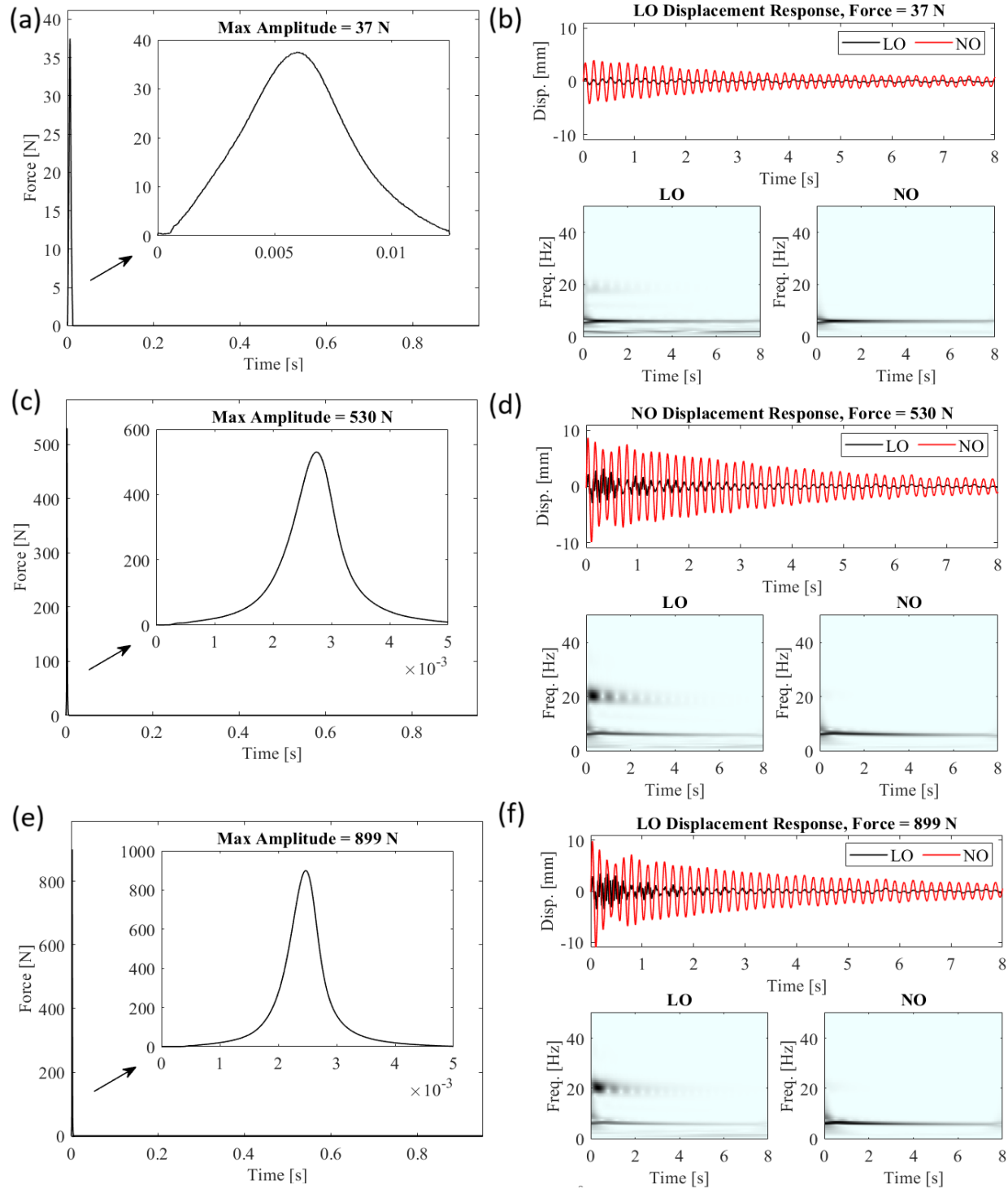

**Fig. 7.** Measured excitation and displacement responses of the LO and NO for impact forces of 37 N, 530 N, and 899 N.

response at later times. We use the time series corresponding to 530 N for training and 37 N and 899 N for validation.

Unlike [67] where the system parameters are identified separately for LO and NO, we use three sets of measurements from the coupled system (LO-NO) to determine all parameters. We use one dataset for identification and two for validation. Additionally, to decrease the computational cost in SIVA, the signals were downsampled to 256 Hz and shifted forward in time to just after the external force. The main frequency content of the system resides in the band $f \in$ [0,50] Hz, such that this downsampling does not remove any important dynamics from the signals. The shifting of time is applied to identify the free-response of the system without needing the external force signal. This ensures that the SIVA methodology can be applied in cases where the external impulsive force cannot be directly measured, such as in shock excitation induced by a blast or acoustic source. The dynamics of the system are modeled using the equations of motion [67]

$$m_{LO}\ddot{x} + b_1\dot{x} + b_2(\dot{x} - \dot{y}) + k_1 x + k_2(x - y) + \alpha(x - y)|x - y|^\beta = 0, \tag{7a}$$

$$m_{NO}\ddot{y} - b_2(\dot{x} - \dot{y}) - k_2(x - y) - \alpha(x - y)|x - y|^\beta = F(t), \tag{7b}$$

where $F(t) = 0$ in the training process of SIVA and $F(t)$ is the externally applied force (i.e., the impact from the modal hammer) when the system is simulated using the identified parameters. The training results of the proposed methodology is presented in Fig. 8(a), showing considerable fluctuations in the early epochs in the losses as the parameter-generator and discriminator compete optimal performance. Convergence occurred after approximately 800 epochs, with all the losses stabilizing.

Once the networks have been trained, each parameter is sampled 1000 times using the parameter-generator, and the resulting mean values are presented in Table 2 as Approach I (this is the same approach used for Approach I as in the simulated system). For the parameters that are in the convergence zone (from epoch 900 onwards) of Fig. 8(b), we obtain the parameter distributions and present them in Figure 8(c). The results suggest that the identified parameters are well-concentrated around their mean value with relatively small standard deviations, indicating a stable and consistent estimation process. Note that the distribution for $\beta$ is not normal, but the individual values are each extremely close to 2. More importantly, the use of random inputs helps prevent the model from overfitting, even with extended training. The mean values of this UQ procedure are presented in Table 2 as Approach II. From the MSE of these two approaches, we can see that the mean values obtained from Approach II provide better results than those from Approach I. To refine the parameters combination, we also present the set of parameters as the result of using Approach III (e.g., obtaining the optimal set of parameters), and the obtained parameters are depicted as Approach III. We calculate MSE as $MSE = \frac{1}{J}\sum_{j=1}^{J}\left[(x_j - \tilde{x}_j)^2 + (y_j - \tilde{y}_j)^2\right]$.

To compare the simulated and measured responses of the LO and NO, we numerically integrate Eq. (7), for the mean values of Table 2 for Approaches II and III, using `ode45` in MATLAB® with zero initial conditions, both relative and absolute and relative tolerances as $10^{-8}$. The measured force is applied up to a time of 0.0062 s using the `interp1` function to interpolate force values between sampled time points in the measured force signal. Figure 9 compares the measured and simulated displacement responses of the LO (Fig. 9 -left panel) and NO (Fig. 9 - right panel). For Approach II, the time series show strong agreement only in early times between the measured and simulated responses, while both the CWT and Fourier spectra show good agreement. As shown in Fig. 9, there is an "optimal" set of parameters (Approach III) from the calculated parameters at each epoch that results in a close match to the number of simulations to perform.

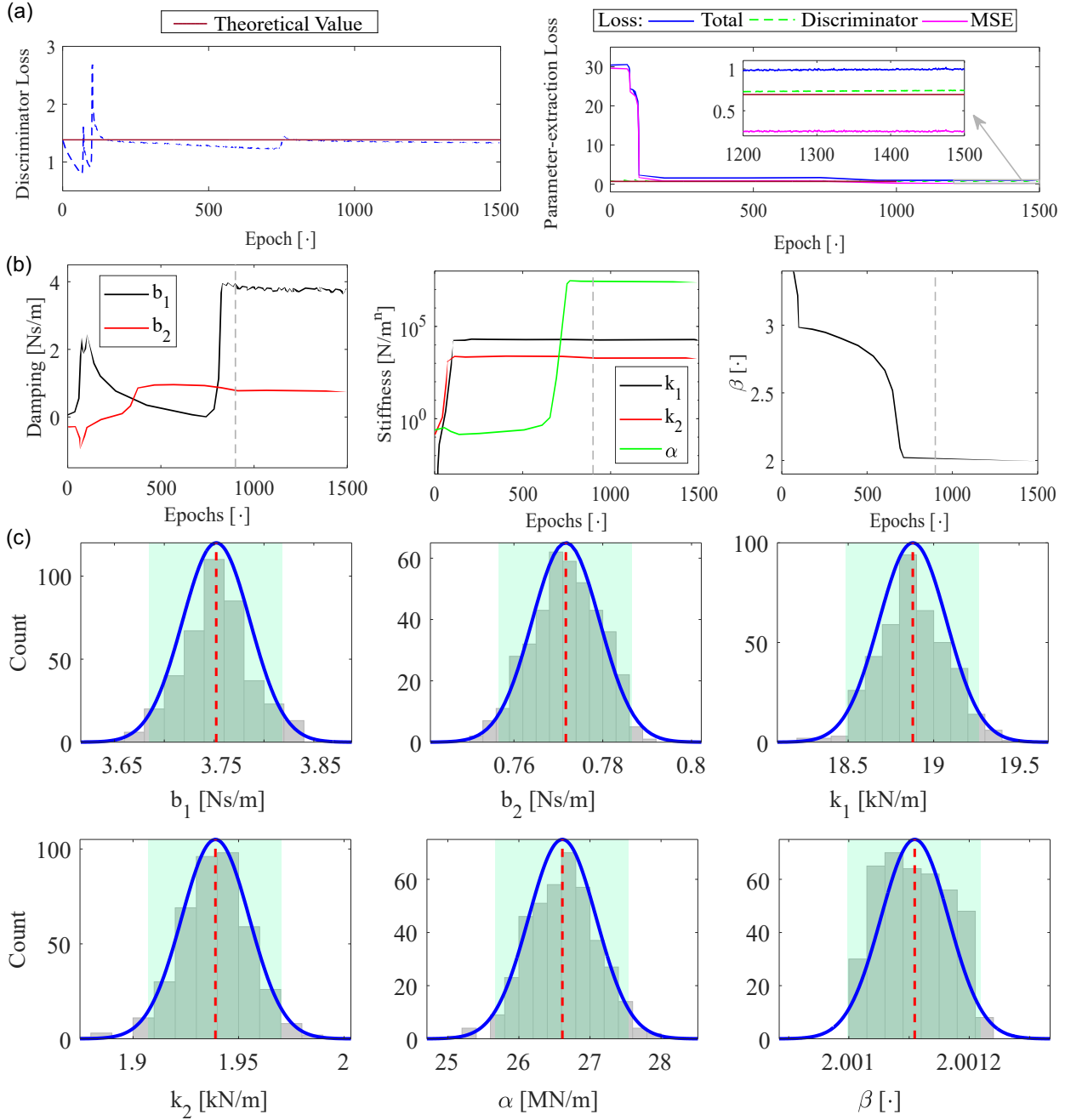

**Fig. 8.** Training dynamics of the proposed SIVA: (a) Losses, (b) Identified parameters distributions of the LO-NO system.

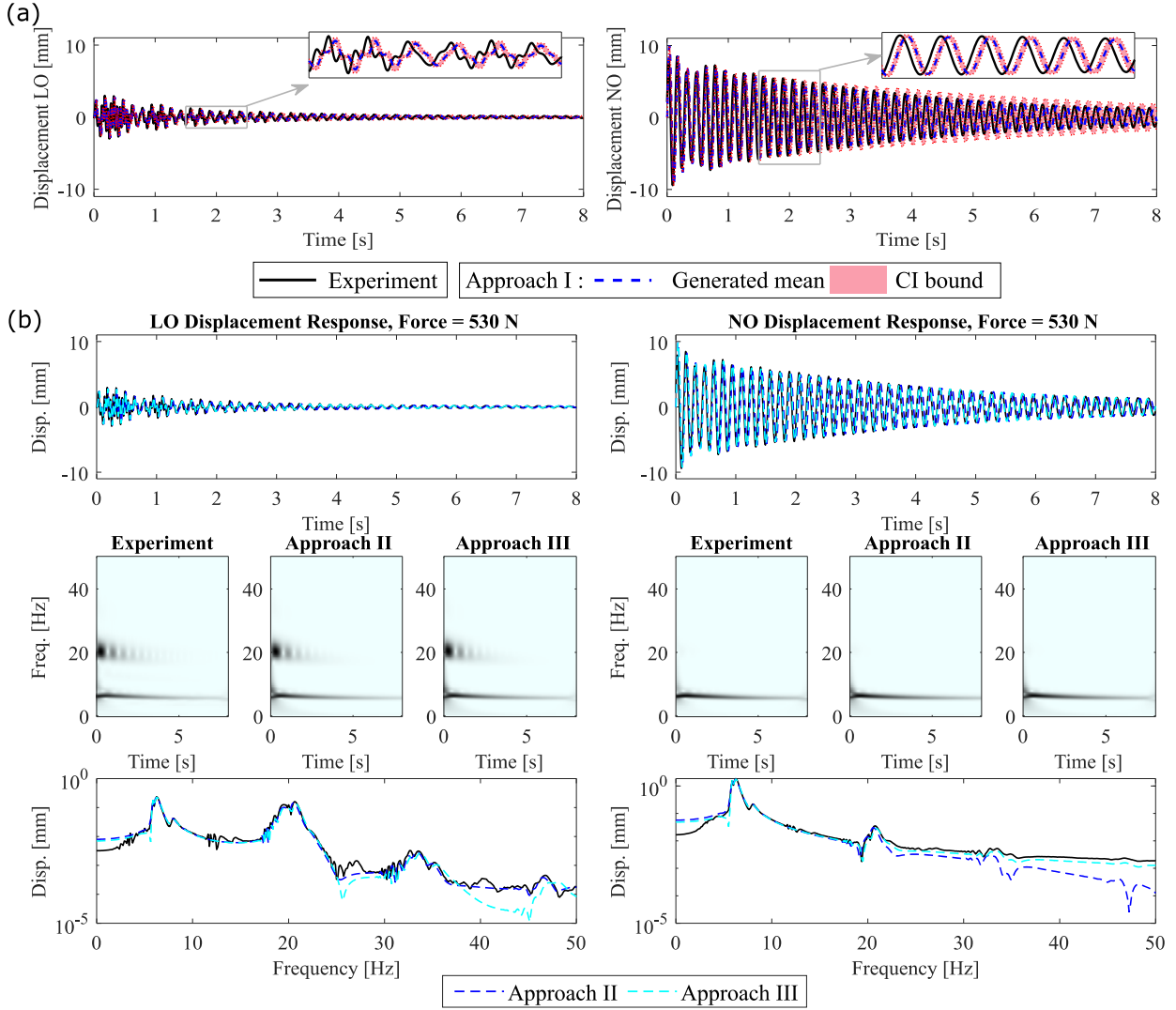

**Fig. 9.** Comparison of the experimentally measured and simulated predicted responses for the (a) LO and (b) NO for the measurement case used in identification by SIVA.

Besides the SIVA method, we utilized the SINDy method with a threshold of 0.1, employing the library of functions $x$, $y$, $\dot{x}$, $\dot{y}$, and $(x-y)^3$. Note that the last term is different from the one presented in Eq. (7) because the standard SINDy only allows standard polynomial terms. However, since we found $\beta = 2$ without much variation, the chosen model for SINDy remains comparable to that identified using SIVA. SINDy resulted in the following parameters values and equation of motion

$$m_{\text{LO}}\ddot{x} + 4.6293\dot{x} - 0.78710\dot{y} + 20535x - 18568y + 2.7260 \times 10^7(x-y)^3 = 0, \quad (8a)$$

$$m_{\text{NO}}\ddot{y} - 1.1905\dot{x} + 0.75080\dot{y} - 1766.6x + 1971.9y - 2.5934 \times 10^7(x-y)^3 = F(t). \quad (8b)$$

Note that Eq. (8) is sorted differently from Eq. (7). This obeys how the library of functions was constructed. We note that implementing constraints in the SINDy method [69] could improve the accuracy of the identified Eq. (7) by constraining the coupling terms to be equal in each equation.

**Table 2.** Comparison of the coefficients of the two-degree-of-freedom model obtained by SIVA.

| Coefficient | Approach I Identified | Approach II Identified | Approach III Identified |
|---|---|---|---|
| $b_1$ [Ns/m] | 3.7270 | 3.7751 | 3.8405 |
| $b_2$ [Ns/m] | 0.75904 | 0.7751 | 0.78563 |
| $k_1$ [N/m] | 18690 | 18883 | 19273 |
| $k_2$ [N/m] | 1920.2 | 1933.8 | 1947.6 |
| $\alpha$ [N/m$^{\beta+1}$] | $2.5815 \times 10^7$ | $2.6805 \times 10^7$ | $2.7641 \times 10^7$ |
| $\beta$ | 2.0010 | 2.0012 | 2.0012 |
| MSE | 10.703 | 5.9052 | 1.3934 |

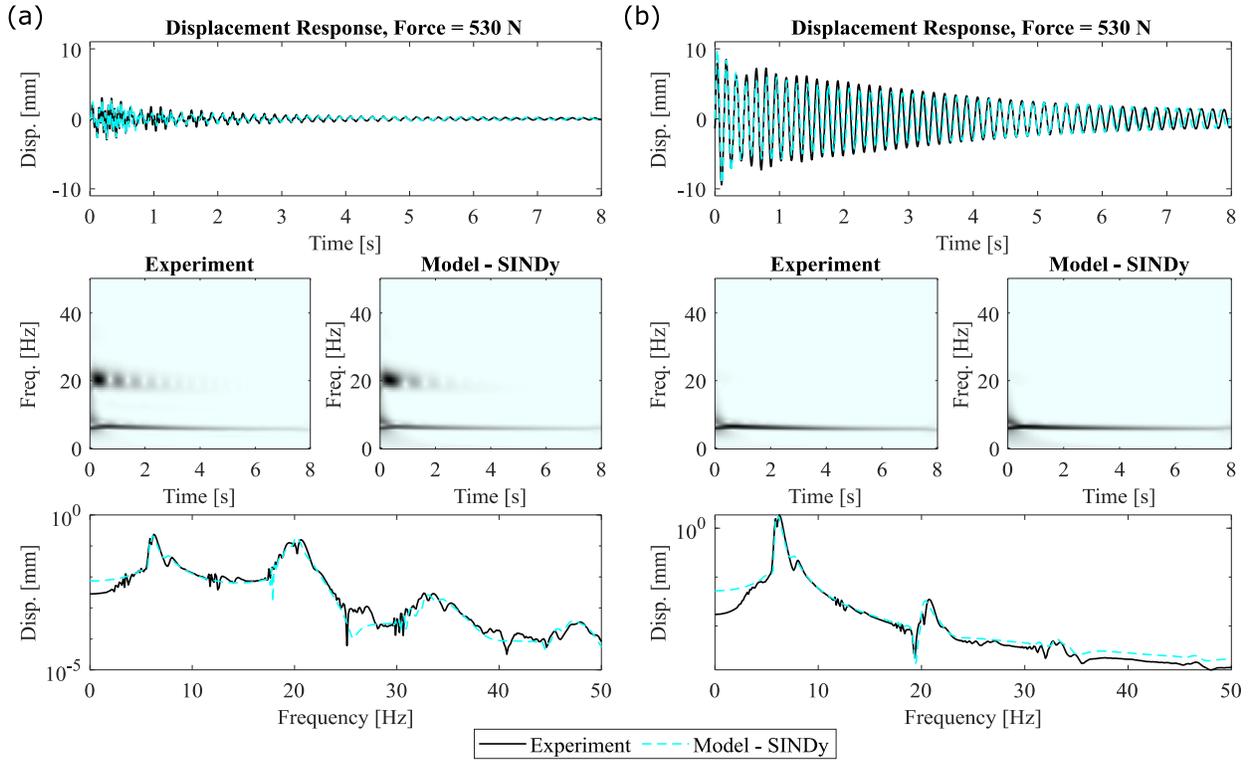

**Fig. 10.** Comparison of the experimentally measured and simulated predicted responses for the (a) LO and (b) NO for the measurement case used in identification by SINDy method.

Figure 10 compares the LO and NO displacement time series, CWT spectra, and Fourier spectra of the measured system with those simulated by integrating Eq. (8), as done in Fig. 9. As can be seen, SINDy provides an acceptable result with a small offset in the linear regime and missing the dynamics for the LO between 1 to 2 seconds. For this system, we do not use the EDDI method, since up to now it only works for single-degree-of-freedom systems.

A further validation of the SINDy and SIVA methods is presented by using the identified parameters in Eq. (8) and Approach III, respectively, to match measurements that were not included in the identification process. To illustrate this, we compare the effect of impacts at 37 N and 899 N in Figs. 11(a) and (b), respectively. For a 37 N impact, the response is weakly nonlinear, and the models for both methods and experimental results align well when comparing the time series and frequency plots. This indicates that the values of damping and stiffness are accurately identified. For the 889 N impact, the response depicts a strongly nonlinear system, evident in both the CWT and Fourier spectra. In the LO plot, the SIVA method captures better the dynamics, as can be seen in the wavelet transform plot. For the NO plot, the simulated responses for both methods reproduce relatively close to the experimental measurements, with a slight offset.

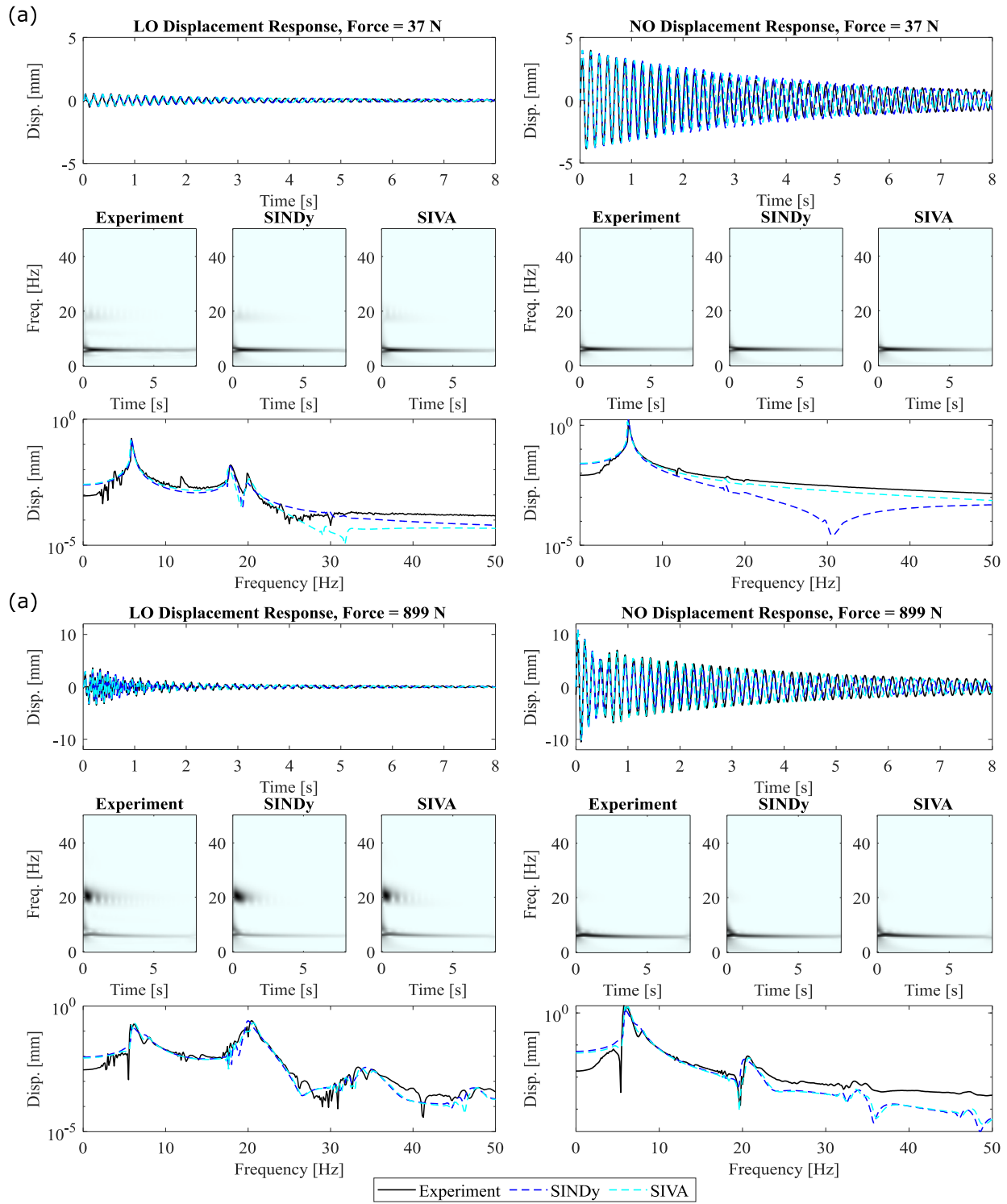

**Fig. 11.** Experimentally measured excitation and displacement responses of the LO-NO for impact amplitudes of (a) 37 N and (b) 899 N.

## 4. Conclusions

We introduced a new parametric system identification method inspired by generative adversarial networks. The so-called structural system identification via validation and adaptation (SIVA) method relies on the measured masses, transient responses, and the mathematical model of the system of interest. The method contains a parameter-generator network that converts random noise into meaningful physical parameter values. The stochastic nature of the input provides regularization during training, helping to reduce overfitting. This is guided by both the mean square error between the real accelerations and those generated by the identified model for the identification dataset, and an adversarial loss that evaluates how well the accelerations generated by the model for validation datasets reproduce their real counterparts. Unlike previous system identification methods, we simultaneously validate the model by using an independent validation dataset. Additionally, we leverage the training procedure to perform uncertainty quantification on the obtained parameters. The proposed approach was applied to both simulated and experimentally measured responses of vibrating oscillators, and the results demonstrated that the SIVA method can perform accurate parametric system identification of strongly nonlinear systems.

While the SIVA shows promising results by identifying the parameters of the tested systems, it relies on the known governing equation of motion. Our studies suggest that the method is relatively robust to the inclusion of additional terms; however, its accuracy may be significantly reduced if key physical terms are missing from the assumed model. This highlights the importance of starting with a reasonably precise mathematical model of the system's dynamics. A future direction is to extend the proposed approach to systems where the governing equations are entirely unknown or only partially known. In the presence of extreme noise, the networks may overfit to patterns caused by the noise rather than learning the underlying physical dynamics, affecting the parameter estimation. Thus, noise is a concern and should be minimized and filtered out as best as possible using bandpass filters. By using the `filtfilt` function in MATLAB, we obtained acceptable results. However, more advanced filtering techniques could be implemented. In this work, the framework has been developed directly from the vanilla GAN, however, GANs face instability in training [70]. In this sense, improvements to GANs could be incorporated into the proposed approach to enhance accuracy and reduce computation time.

## Acknowledgments

This research was supported by the Air Force Office of Scientific Research Young Investigator Program under grant number FA9550-22-1-0295.
## References

1. Worden, K., Tomlinson, G.R.: Nonlinearity in Structural Dynamics. Institute of Physics Publishing, Bristol and Philadelphia (2001)
2. Malekloo, A., Ozer, E., AlHamaydeh, M., Girolami, M.: Machine learning and structural health monitoring overview with emerging technology and high-dimensional data source highlights. Struct. Heal. Monit. 21, 1906–1955 (2022). https://doi.org/10.1177/14759217211036880
3. Kerschen, G., Worden, K., Vakakis, A.F., Golinval, J.C.: Past, present and future of nonlinear system identification in structural dynamics. Mech. Syst. Signal Process. 20, 505–592 (2006). https://doi.org/10.1016/j.ymssp.2005.04.008
4. Noël, J.P., Kerschen, G.: Nonlinear system identification in structural dynamics: 10 more years of progress. Mech. Syst. Signal Process. 83, 2–35 (2017). https://doi.org/10.1016/j.ymssp.2016.07.020
5. Bunce, A., Brennan, D.S., Ferguson, A., O'Higgins, C., Taylor, S., Cross, E.J., Worden, K., Brownjohn, J., Hester, D.: On population-based structural health monitoring for bridges: Comparing similarity metrics and dynamic responses between sets of bridges. Mech. Syst. Signal Process. 216, 111501 (2024). https://doi.org/10.1016/j.ymssp.2024.111501
6. Alkhatib, R., Golnaraghi, M.F.: Active structural vibration control: A review. Shock Vib. Dig. 35, 367 (2003). https://doi.org/10.1177/05831024030355002
7. Bagha, A.K., Modak, S. V.: Active structural-acoustic control of interior noise in a vibro-acoustic cavity incorporating system identification. J. Low Freq. Noise Vib. Act. Control. 36, 261–276 (2017). https://doi.org/10.1177/0263092317719636
8. Ewins, D.J.: Modal Testing: Theory, Practice, and Application. Research Studies Press (2000)
9. Keesman, K.J.: Sytem Identification, An Introduction. Springer London, New York (2011)
10. Guan, W., Dong, L.L., Zhou, J.M., Han, Y., Zhou, J.: Data-driven methods for operational modal parameters identification: A comparison and application. Meas. J. Int. Meas. Confed. 132, 238–251 (2019). https://doi.org/10.1016/j.measurement.2018.09.052
11. Sadhu, A., Narasimhan, S., Antoni, J.: A review of output-only structural mode identification literature employing blind source separation methods. Mech. Syst. Signal Process. 94, 415–431 (2017). https://doi.org/10.1016/j.ymssp.2017.03.001
12. Lai, Z., Nagarajaiah, S.: Sparse structural system identification method for nonlinear dynamic systems with hysteresis/inelastic behavior. Mech. Syst. Signal Process. 117, 813–842 (2019). https://doi.org/10.1016/j.ymssp.2018.08.033
13. Peng, Z.K., Lang, Z.Q., Wolters, C., Billings, S.A., Worden, K.: Feasibility study of structural damage detection using NARMAX modelling and Nonlinear Output Frequency Response Function based analysis. Mech. Syst. Signal Process. 25, 1045–1061 (2011). https://doi.org/10.1016/j.ymssp.2010.09.014
14. Leontaritis, I.J., Billings, S.A.: Input-output parametric models for non-linear systems Part
21

I: Deterministic non-linear systems. Int. J. Control. 41, 303–328 (1985). https://doi.org/10.1080/0020718508961129
15. Chatzi, E.N., Smyth, A.W., Masri, S.F.: Experimental application of on-line parametric identification for nonlinear hysteretic systems with model uncertainty. Struct. Saf. 32, 326–337 (2010). https://doi.org/10.1016/j.strusafe.2010.03.008
16. Ben Abdessalem, A., Dervilis, N., Wagg, D., Worden, K.: Model selection and parameter estimation in structural dynamics using approximate Bayesian computation. Mech. Syst. Signal Process. 99, 306–325 (2018). https://doi.org/10.1016/j.ymssp.2017.06.017
17. Ben Abdessalem, A., Dervilis, N., Wagg, D., Worden, K.: Model selection and parameter estimation of dynamical systems using a novel variant of approximate Bayesian computation. Mech. Syst. Signal Process. 122, 364–386 (2019). https://doi.org/10.1016/j.ymssp.2018.12.048
18. Scheel, M., Kleyman, G., Tatar, A., Brake, M.R.W., Peter, S., Noël, J.P., Allen, M.S., Krack, M.: Experimental assessment of polynomial nonlinear state-space and nonlinear-mode models for near-resonant vibrations. Mech. Syst. Signal Process. 143, 106796 (2020). https://doi.org/10.1016/j.ymssp.2020.106796
19. Moore, K.J., Bunyan, J., Tawfick, S., Gendelman, O. V., Li, S., Leamy, M., Vakakis, A.F.: Nonreciprocity in the dynamics of coupled oscillators with nonlinearity, asymmetry, and scale hierarchy. Phys. Rev. E. 97, 1–11 (2018). https://doi.org/10.1103/PhysRevE.97.012219
20. Bunyan, J., Moore, K.J., Mojahed, A., Fronk, M.D., Leamy, M., Tawfick, S., Vakakis, A.F.: Acoustic nonreciprocity in a lattice incorporating nonlinearity, asymmetry, and internal scale hierarchy: Experimental study. Phys. Rev. E. 97, 1–13 (2018). https://doi.org/10.1103/PhysRevE.97.052211
21. Lejarza, F., Baldea, M.: Data-driven discovery of the governing equations of dynamical systems via moving horizon optimization. Sci. Rep. 12, 1–15 (2022). https://doi.org/10.1038/s41598-022-13644-w
22. Zhang, R., Liu, Y., Sun, H.: Physics-informed multi-LSTM networks for metamodeling of nonlinear structures. Comput. Methods Appl. Mech. Eng. 369, 113226 (2020). https://doi.org/10.1016/j.cma.2020.113226
23. Ding, Z., Yu, Y., Xia, Y.: Nonlinear hysteretic parameter identification using an attention-based long short-term memory network and principal component analysis. Nonlinear Dyn. 111, 4559–4576 (2023). https://doi.org/10.1007/s11071-022-08095-x
24. Zhai, W., Tao, D., Bao, Y.: Parameter estimation and modeling of nonlinear dynamical systems based on Runge–Kutta physics-informed neural network. Nonlinear Dyn. 111, 21117–21130 (2023). https://doi.org/10.1007/s11071-023-08933-6
25. Yu, Y., Liu, Y.: Physics-guided generative adversarial network for probabilistic structural system identification. Expert Syst. Appl. 239, 122339 (2024). https://doi.org/10.1016/j.eswa.2023.122339
26. Liu, T., Meidani, H.: Physics-Informed Neural Networks for System Identification of Structural Systems with a Multiphysics Damping Model. J. Eng. Mech. 149, 1–12 (2023). https://doi.org/10.1061/jenmdt.emeng-7060
27. Masri, S.F., Caughey, T.K.: A nonparametric identification technique for nonlinear dynamic problems. J. Appl. Mech. Trans. ASME. 46, 433–447 (1979). https://doi.org/10.1115/1.3424568
28. Noël, J.P., Renson, L., Kerschen, G.: Complex dynamics of a nonlinear aerospace

structure: Experimental identification and modal interactions. J. Sound Vib. 333, 2588–2607 (2014). https://doi.org/10.1016/j.jsv.2014.01.024
29. Gray, G.J., Murray-Smith, D.J., Li, Y., Sharman, K.C., Weinbrenner, T.: Nonlinear model structure identification using genetic programming. Control Eng. Pract. 6, 1341–1352 (1998). https://doi.org/10.1016/S0967-0661(98)00087-2
30. Bolourchi, A., Masri, S.F., Aldraihem, O.J.: Studies into computational intelligence and evolutionary approaches for model-free identification of hysteretic systems. Comput. Civ. Infrastruct. Eng. 30, 330–346 (2015). https://doi.org/10.1111/mice.12126
31. Im, J., Rizzo, C.B., de Barros, F.P.J., Masri, S.F.: Application of genetic programming for model-free identification of nonlinear multi-physics systems. Nonlinear Dyn. 104, 1781–1800 (2021). https://doi.org/10.1007/s11071-021-06335-0
32. Masri, S.F., Chassiakos, A.G., Caughey, T.K.: Identification of nonlinear dynamic systems using neural networks. J. Appl. Mech. 60, 123–133 (1993). https://doi.org/https://doi.org/10.1115/1.2900734
33. Adeli, H., Jiang, X.: Dynamic Fuzzy Wavelet Neural Network for Structural System Identification. Intell. Infrastruct. 132, 102–111 (2006). https://doi.org/doi.org/10.1061/(ASCE)0733-9445(2006)132:1(102)
34. Zhou, J.M., Dong, L., Guan, W., Yan, J.: Impact load identification of nonlinear structures using deep Recurrent Neural Network. Mech. Syst. Signal Process. 133, 106292 (2019). https://doi.org/10.1016/j.ymssp.2019.106292
35. Wu, R.-T., Jahanshahi, M.R.: Deep Convolutional Neural Network for Structural Dynamic Response Estimation and System Identification. J. Eng. Mech. 145, (2019). https://doi.org/10.1061/(asce)em.1943-7889.0001556
36. Chen, Z., Liu, Y., Sun, H.: Symbolic Deep Learning for Structural System Identification. J. Struct. Eng. 148, 1–14 (2022). https://doi.org/10.1061/(asce)st.1943-541x.0003405
37. Lee, Y.S., Vakakis, A.F., McFarland D.M., Bergman, L.A.: A global–local approach to nonlinear system identification: A review. Struct. Control Heal. Monit. 17, 742–760 (2010). https://doi.org/10.1002/stc
38. Vakakis, A.F., Bergman, L.A., McFarland, D.M., Lee, Y.S., Kurt, M.: Current efforts towards a non-linear system identification methodology of broad applicability. Proc. Inst. Mech. Eng. Part C J. Mech. Eng. Sci. 225, 2497–2515 (2011). https://doi.org/10.1177/0954406211417217
39. Kwarta, M., Allen, M.S.: NIXO-Based identification of the dominant terms in a nonlinear equation of motion of structures with geometric nonlinearity. J. Sound Vib. 568, (2024). https://doi.org/10.1016/j.jsv.2023.117900
40. Schmidt, M., Lipson, H.: Distilling free-form natural laws from experimental data. Science (80-. ). 324, 81–85 (2009). https://doi.org/10.1126/science.1165893
41. Liu, Z., Tegmark, M.: Machine Learning Conservation Laws from Trajectories. Phys. Rev. Lett. 126, 180604 (2021). https://doi.org/10.1103/PhysRevLett.126.180604
42. Allen, M.S., Sumali, H., Epp, D.S.: Piecewise-linear restoring force surfaces for semi-nonparametric identification of nonlinear systems. Nonlinear Dyn. 54, 123–135 (2008). https://doi.org/10.1007/s11071-007-9254-x
43. Brunton, S.L., Proctor, J.L., Kutz, J.N.: Discovering governing equations from data by sparse identification of nonlinear dynamical systems. Proc. Natl. Acad. Sci. U. S. A. 113, 3932–3937 (2016). https://doi.org/10.1073/pnas.1517384113
44. Fuentes, R., Nayek, R., Gardner, P., Dervilis, N., Rogers, T., Worden, K., Cross, E.J.:

Equation discovery for nonlinear dynamical systems: A Bayesian viewpoint. Mech. Syst. Signal Process. 154, 107528 (2021). https://doi.org/10.1016/j.ymssp.2020.107528
45. Moore, K.J.: Characteristic nonlinear system identification: A data-driven approach for local nonlinear attachments. Mech. Syst. Signal Process. 131, 335–347 (2019). https://doi.org/10.1016/j.ymssp.2019.05.066
46. Najera-Flores, D.A., Todd, M.D.: A structure-preserving neural differential operator with embedded Hamiltonian constraints for modeling structural dynamics. Comput. Mech. 72, 241–252 (2023). https://doi.org/10.1007/s00466-023-02288-w
47. López, C., Singh, A., Naranjo, Á., Moore, K.J.: A data-driven, energy-based approach for identifying equations of motion in vibrating structures directly from measurements. Mech. Syst. Signal Process. 225, 112341 (2025). https://doi.org/10.1016/j.ymssp.2025.112341
48. López, C., Moore, K.J.: Energy-based dual-phase dynamics identification of clearance nonlinearities. Nonlinear Dyn. 113, 17933–17948 (2025). https://doi.org/10.1007/s11071-025-11098-z
49. Quaranta, G., Lacarbonara, W., Masri, S.F.: A review on computational intelligence for identification of nonlinear dynamical systems. Nonlinear Dyn. 99, 1709–1761 (2020). https://doi.org/10.1007/s11071-019-05430-7
50. Zaparoli Cunha, B., Droz, C., Zine, A.M., Foulard, S., Ichchou, M.: A review of machine learning methods applied to structural dynamics and vibroacoustic. Mech. Syst. Signal Process. 200, 110535 (2023). https://doi.org/10.1016/j.ymssp.2023.110535
51. Goodfellow, I.J., Pouget-Abadie, J., Mirza, M., Xu, B., Warde-Farley, D., Ozair, S., Courville, A., Bengio, Y.: Generative adversarial nets. Adv. Neural Inf. Process. Syst. 27, 2672–2680 (2014)
52. Goodfellow, I., Pouget-Abadie, J., Mirza, M., Xu, B., Warde-Farley, D., Ozair, S., Courville, A., Bengio, Y.: Generative adversarial networks. Commun. ACM. 63, 139–144 (2020). https://doi.org/10.1145/3422622
53. Yu, L., Zhang, W., Wang, J., Yu, Y.: SeqGAN: Sequence generative adversarial nets with policy gradient. Proc. AAAI Conf. Artif. Intell. 31, (2017). https://doi.org/10.1609/aaai.v31i1.10804
54. de Rosa, G.H., Papa, J.P.: A survey on text generation using generative adversarial networks. Pattern Recognit. 119, 108098 (2021). https://doi.org/10.1016/j.patcog.2021.108098
55. Shao, S., Wang, P., Yan, R.: Generative adversarial networks for data augmentation in machine fault diagnosis. Comput. Ind. 106, 85–93 (2019). https://doi.org/10.1016/j.compind.2019.01.001
56. Akcay, S., Atapour-Abarghouei, A., Breckon, T.P.: Ganomaly: Semi-supervised Anomaly Detection via Adversarial Training. In: Asian conference on computer vision. pp. 622–637. Springer, Cham, Australia, 2-6 December (2018)
57. Zhong, C., Zhang, J., Lu, X., Zhang, K., Liu, J., Hu, K., Chen, J., Lin, X.: Deep Generative Model for Inverse Design of High-Temperature Superconductor Compositions with Predicted Tc > 77 K. ACS Appl. Mater. Interfaces. 15, 30029–30038 (2023). https://doi.org/10.1021/acsami.3c00593
58. Yang, L., Zhang, D., Karniadakis, G.E.: Physics-informed generative adversarial networks for stochastic differential equations. SIAM J. Sci. Comput. 42, A292-317 (2020)
59. Pan, T., Chen, J., Zhang, T., Liu, S., He, S., Lv, H.: Generative adversarial network in mechanical fault diagnosis under small sample: A systematic review on applications and

future perspectives. ISA Trans. 128, 1–10 (2022). https://doi.org/10.1016/j.isatra.2021.11.040
60. Yuan, Z.Q., Xin, Y., Wang, Z.C., Ding, Y.J., Wang, J., Wang, D.H.: Structural Nonlinear Model Updating Based on an Improved Generative Adversarial Network. Struct. Control Heal. Monit. 2023, 1–21 (2023). https://doi.org/10.1155/2023/9278389
61. Rostamijavanani, A., Li, S., Yang, Y.: Data-Driven Modeling of Parameterized Nonlinear Dynamical Systems with a Dynamics-Embedded Conditional Generative Adversarial Network. J. Eng. Mech. 149, 1–13 (2023). https://doi.org/10.1061/jenmdt.emeng-7038
62. Tsialiamanis, G., Champneys, M.D., Dervilis, N., Wagg, D.J., Worden, K.: On the application of generative adversarial networks for nonlinear modal analysis. Mech. Syst. Signal Process. 166, 108473 (2022). https://doi.org/10.1016/j.ymssp.2021.108473
63. Moore, K.J., Kurt, M., Eriten, M., McFarland, D.M., Bergman, L.A., Vakakis, A.F.: Direct detection of nonlinear modal interactions from time series measurements. Mech. Syst. Signal Process. 125, 311–329 (2019). https://doi.org/10.1016/j.ymssp.2017.09.010
64. Goodfellow, I.: NIPS 2016 Tutorial: Generative Adversarial Networks. (2016)
65. Paul S. Addison: The Illustrated wavelet Transform Handbook Introductory Theory and Applications in Science, Engineering, Medicine and Finance. , Boca Raton (2017)
66. Moore, K.J., Kurt, M., Eriten, M., McFarland, D.M., Bergman, L.A., Vakakis, A.F.: Wavelet-bounded empirical mode decomposition for measured time series analysis. Mech. Syst. Signal Process. 99, 14–29 (2018). https://doi.org/10.1016/j.ymssp.2017.06.005
67. Wang, C., Moore, K.J.: On nonlinear energy flows in nonlinearly coupled oscillators with equal mass. Nonlinear Dyn. 103, 343–366 (2021). https://doi.org/10.1007/s11071-020-06120-5
68. Wang, C., González, G.Y., Wittich, C., Moore, K.J.: Energy isolation in a multi-floor nonlinear structure under harmonic excitation. Nonlinear Dyn. 110, 2049–2077 (2022). https://doi.org/10.1007/s11071-022-07744-5
69. Loiseau, J.C., Brunton, S.L.: Constrained sparse Galerkin regression. J. Fluid Mech. 838, 42–67 (2018). https://doi.org/10.1017/jfm.2017.823
70. Arjovsky, M., Chintala, S., Bottou, L.: Wasserstein GAN. In: International conference on machine learning. pp. 214–223. PMLR (2017)